\documentclass[review]{elsarticle}
\RequirePackage{amssymb}
\usepackage{bm}
\RequirePackage{amsmath}
\usepackage{amsthm}
\usepackage{color}
\usepackage{lineno,hyperref}
\modulolinenumbers[5]
\usepackage{geometry}
\usepackage{lineno,hyperref}
\modulolinenumbers[5]

\begin{document}
\newtheorem{theorem}{Theorem}[section]
\newtheorem{definition}{Definition}[section]
\newtheorem{lemma}{Lemma}[section]
\newtheorem{problem}{Problem}[section]
\newtheorem{claim}{Claim}[section]
\newtheorem{corollary}{Corollary}[section]
\newtheorem{proposition}{Proposition}[section]
\newtheorem{conjecture}{Conjecture}[section]
\newtheorem{observation}{Observation}[section]
 \abovedisplayskip 6pt plus 2pt minus 2pt \belowdisplayskip 6pt
 plus 2pt minus 2pt
 %%%%%%%%%%%%%%%%
\def\vsp{\vspace{1mm}}
\def\th#1{\vspace{1mm}\noindent{\bf #1}\quad}
\def\proof{\vspace{1mm}\noindent{\it Proof}\quad}
\def\no{\nonumber}
\newenvironment{prof}[1][Proof]{\noindent\textit{#1}\quad }
{\hfill $\Box$\vspace{0.7mm}}
\begin{frontmatter}

\title{Anti-Ramsey numbers of loose paths and cycles in uniform hypergraphs}

\author[ams,ucas]{Tong Li}
\ead{litong@amss.ac.cn}

\author[nuaa,klmm]{Yucong Tang}
\ead{tangyucong@nuaa.edu.cn}

\author[sdu]{Guanghui Wang}
\ead{ghwang@sdu.edu.cn}

\author[ams,ucas]{Guiying Yan\corref{cor1}}
\ead{yangy@amss.ac.cn}
%\fntext[NSF1]{Research partly supported by  National Natural Science Foundation of China (Grant No. 11631014)}

\cortext[cor1]{Corresponding author}
\address[ams]
{Academy of Mathematics and Systems Science, Chinese Academy of Sciences, Beijing $100190$, P. R. China}

\address[nuaa]
{Department of Mathematics, Nanjing University of Aeronautics and Astronautics, Nanjing $211106$, P. R. China}

\address[klmm]
{Key Laboratory of Mathematical Modeling and High Performance Computing of Air Vehicles (NUAA), MIIT, Nanjing $211106$, P. R. China}

\address[sdu]
{School of Mathematics, Shandong University, Jinan $250100$, P. R. China}

\address[ucas]
{School of Mathematical Sciences, University of Chinese Academy of Sciences, Beijing $100049$, P. R. China}

\begin{abstract}

For a fixed family of $r$-uniform hypergraphs $\mathcal{F}$,  the anti-Ramsey number of $\mathcal{F}$, denoted by $ ar(n,r,\mathcal{F})$,
is the minimum number $c$ of colors such that for any edge-coloring of the complete $r$-uniform hypergraph on $n$ vertices with at least $c$ colors, there is a rainbow copy of some hypergraph in $\mathcal{F}$. Here, a rainbow hypergraph is an edge-colored hypergraph with all edges colored differently. Let $\mathcal{P}_k$ and $\mathcal{C}_k$ be the families of loose paths and loose cycles with $k$ edges in an $r$-uniform hypergraph, respectively. In this paper, we determine the exact values of $ ar(n,r,\mathcal{P}_k)$ and $ ar(n,r,\mathcal{C}_k)$ for all $k\geq 4$ and $r\geq 3$.

\end{abstract}


\end{frontmatter}

%\linenumbers

\section{Introduction}

An \emph{$r$-uniform hypergraph} (or $r$-graph, for simplicity) is a pair $H=(V, E)$, where $V=V(H)$ is a finite set of vertices, and $E=E(H)\subseteq {{V}\choose {r}}$ is a family of $r$-element subsets of $V$ called edges. Let $K_{n}^{r}$ be the complete $r$-graph on $n$ vertices. For a fixed family of $r$-graphs $\mathcal{F}$,  the\emph{ anti-Ramsey number} of $\mathcal{F}$, denoted by $ar(n,r,\mathcal{F})$, is the minimum number $c$ of colors such that for any edge-coloring of $K_{n}^{r}$ with at least $c$ colors, there is a rainbow copy of some hypergraph in $\mathcal{F}$. Here, a \emph{rainbow} hypergraph is an edge-colored hypergraph with all edges colored differently.

Anti-Ramsey numbers were introduced by Erd\H{o}s, Simonovits and S$\acute{\rm o}$s  \cite{Erdos&Simonovits&Sos1975} in 1973. They found that the anti-Ramsey numbers are closely related to the Tur\'{a}n numbers. The \emph{Tur\'{a}n number} of a family of $r$-graphs $\mathcal{F}$, denoted by $ex(n,r,\mathcal{F})$, is the maximum number of edges in an $r$-graph on $n$ vertices that does not contain any $r$-graph in $\mathcal{F}$ as a subgraph. For a fixed $r$-graph $F$, there is a natural lower bound of $ar(n,r,F)$ in terms of  Tur\'{a}n number as follow.
\begin{observation}\label{trivial}
$$ar(n,r,F)\geq ex(n,r,\{F-e:e\in F\})+2.$$
\end{observation}
This trivial lower bound is easily obtained by coloring a rainbow Tur\'{a}n extremal $r$-graph for $\{F-e:e\in F\}$ in $K_n^r$, and the remaining edges with an additional color.

For $2$-graphs, various results about the anti-Ramsey numbers have been obtained, including those of cliques, paths, cycles, matchings and so on. We refer the readers to the survey \cite{Fujita&Magnant&Ozeki2010}.

For $r$-graphs with $r\geq 3$, not too many results were obtained. \"{O}zkahya and Young \cite{Ozkahya&Young2013} initiated the study of anti-Ramsey number for matchings in $r$-graphs, where a matching $M_k$ is a collection of pairwise disjoint $k$ edges.  They conjectured that
\begin{align*}ar(n,r,M_k)=\left\{\begin{array}{ll}
ex(n, r,M_{k-1}) + 2 , & if\ k< c_r,n=kr,\\
ex(n, r,M_{k-1}) + r+1 , & if\ k\geq c_r,n=kr,\\
ex(n, r,M_{k-1}) + 2, & if\ n>kr,
\end{array}\right.\end{align*}
 where $c_r$ is a constant depending on $r$.  They proved that the conjecture holds for $k=2,3$ and  sufficiently large $n$. Later, using shifting method, Frankl and Kupavskii \cite{Frankl&Kupavskii2019} proved that the conjecture is true for $ n \geq rk + (r - 1)(k - 1)$ and $k \geq 3$. However, this conjecture was proved to be false for $n=kr$ by Guo, Lu and Peng \cite{Guo&Lu&Peng2022} very recently.

Let $[n]=\{1, 2, \ldots, n\}$. A \emph{loose cycle of length} $k$ in an $r$-graph is a collection of $k$ edges $\{e_1,\ldots, e_k\}$ such that $e_i\cap e_{i+1}\neq \emptyset$ for $i\in [k]$ (where $e_{k+1}=e_1$), and $e_i\cap e_j=\emptyset$ for $i-j\not\equiv \pm 1(\text{mod } k)$. Denoted by $\mathcal{C}_{k}$ all loose cycles of length $k$. A \emph{loose path of length} $k$ in an $r$-graph is a collection of $k$ edges $\{e_1,\ldots, e_k\}$ such that $e_i\cap e_{i+1}\neq \emptyset$ for $i\in [k-1]$, and $e_i\cap e_j=\emptyset$ for $i-j\not\equiv \pm 1(\text{mod } k)$. Denoted by $\mathcal{P}_{k}$ all loose paths of length $k$. Gu, Li and Shi \cite{Gu&Li&Shi2020} determined the exact anti-Ramsey numbers for $\mathcal{C}_{k}$ and $\mathcal{P}_{k}$ in the following cases.

\begin{theorem}[\cite{Gu&Li&Shi2020}]
Let $n,r,k,t$ be integers such that $n$ is sufficiently large. Then
\begin{align*}ar(n,r,\mathcal{C}_k)=ar(n,r,\mathcal{P}_k)=\left\{\begin{array}{ll}
{n\choose r}-{{n-t+1}\choose r}+2 , & if\ k=2t\geq 8\ and\ r\geq 4,\\
{n\choose r}-{{n-t+1}\choose r}+3 , & if\ k=2t+1\geq 11\ and\ r\geq k+3.
\end{array}\right.\end{align*}
\end{theorem}

In this paper, we extend their results and determine the exact values of $ar(n,r,\mathcal{P}_k)$ and $ar(n,r,\mathcal{C}_k)$ for all $k\geq 4$ and $r\geq 3$. With the case $k=3$ and $r\geq 3$ having already been considered before, see \cite{Gu&Li&Shi2020,Tang&Li&Yan2023}, the exact values of $ar(n,r,\mathcal{P}_k)$ and $ar(n,r,\mathcal{C}_k)$ for all $k\geq 3$ and $r\geq 3$ are obtained.

\begin{theorem}\label{main theorem}
Let $n,r,k,t$ be integers such that $n$ is sufficiently large and $r\geq 3$. Then
\begin{align*}ar(n,r,\mathcal{C}_k)=ar(n,r,\mathcal{P}_k)=\left\{\begin{array}{ll}
{n\choose r}-{{n-t+1}\choose r}+2 , & if\ k=2t\geq4,\\
{n\choose r}-{{n-t+1}\choose r}+3 , & if\ k=2t+1\geq5.
\end{array}\right.\end{align*}
\end{theorem}

To express the constant hierarchies, we write $x\ll y$ to mean that for any $y\in (0,1]$, there exists an $x_0\in (0,1)$, where $x_0$ is often regarded as a positive function of $y$, such that for all $x\leq x_0$, the subsequent statements hold. Hierarchies with more constants are defined in a similar way and to be read from the right to the left.

The rest of this paper is organized as follows. In Section 2, we give necessary definitions and lemmas for the proofs. In Section 3, we prove the case when $k$ is even. In section 4, we prove the case when $k$ is odd.

\section{Preliminaries}

 A \emph{linear cycle of length} $k$ in an $r$-graph, denoted by  $C_k^+$, is a collection of $k$ edges $\{e_1,\ldots, e_k\}$ such that $|e_i\cap e_{i+1}|=1$ for $i\in [k]$ (where $e_{k+1}=e_1$), and $e_i\cap e_j=\emptyset$ for $i-j\not\equiv \pm 1(\text{mod } k) $. A \emph{linear path of length} $k$ in an $r$-graph, denoted by $P_k^+$, is a collection of $k$ edges $\{e_1,\ldots, e_k\}$ such that $|e_i\cap e_{i+1}|=1$ for $i\in [k-1]$, and $e_i\cap e_j=\emptyset$ for $i-j\not\equiv \pm 1(\text{mod } k) $. Notice that a linear path or cycle is also a loose path or cycle.

In the following, we list some lemmas which will be used in Section 3 and Section 4, including anti-Ramsey numbers of linear paths or cycles and  short loose paths, Tur\'{a}n numbers of linear paths or cycles and loose paths or cycles, and stability lemma.
\begin{lemma}[\cite{Gu&Li&Shi2020, Tang&Li&Yan2022}]\label{linear}
Let $n,r,k,t$ be integers such that $1/n\ll r,k$ and $r\geq 3$. Then
\begin{align*}ar(n,r,C_k^+)=ar(n,r,P_k^+)=\left\{\begin{array}{ll}
{n\choose r}-{{n-t+1}\choose r}+2 , & if\ k=2t\geq4,\\
{n\choose r}-{{n-t+1}\choose r}+{{n-t-1}\choose r-2}+2 , & if\ k=2t+1\geq5.
\end{array}\right.\end{align*}
\end{lemma}

\begin{lemma}[\cite{Gu&Li&Shi2020}]\label{arrP3}
For $n\geq 3r-4$, $ar(n,r,\mathcal{P}_{2})=2$. For $n\geq 4r-3$, $ar(n,r,\mathcal{P}_{3})=3$.
\end{lemma}

\begin{lemma}[\cite{Erdos&Gallai1959}]\label{Turan2Pk}(\textbf{Erd\H{o}s-Gallai Theorem})
For any $k\geq 2$ and  $n\geq 1$,
$ex(n,2,P_k)\leq \frac{k-1}{2}n.
$
\end{lemma}

The determination of $ex(n,r,P_k^+)$ is nontrivial even for $k=2$. We refer the readers to \cite{Frankl1977,Keevash&Mubayi&Wilson2006,Furedi&Jiang&Seiver2014,Jackowska&Polcyn&Rucinski2016} for more results and details. In particular, Kostochka, Mubayi and Verstra\"{e}te \cite{Kostochka&Mubayi&Verstraete2015} considered $ex(n,r,P_k^+)$ for $r\geq 3$, $k\geq 4$ and sufficiently large $n$. F\"{u}redi and Jiang \cite{Furedi&Jiang2014} and Kostochka, Mubayi and Verstra\"{e}te \cite{Kostochka&Mubayi&Verstraete2015} determined $ex(n,r,C_k^+)$ for all $k\geq 3$ and $r\geq 3$ and sufficiently large $n$.

\begin{lemma}[\cite{Furedi&Jiang2014,Kostochka&Mubayi&Verstraete2015}]\label{turanPkCk} Let $\delta,$ $\varepsilon>0$ and $r, k, t, n$ be positive integers such that $r\geq 3,$ $k\geq 4$, $t=\lfloor \frac{k-1}{2}\rfloor$ and $1/n\ll \delta\ll \varepsilon$. Then
\begin{align*}ex(n,r,P_k^+)={{n}\choose {r}}-{{n-t}\choose r}+\left\{\begin{array}{ll}
0, & \mbox{if $k$ is odd,}\\ {{n-t-2}\choose {r-2}}, & \mbox{if $k$ is even}.
\end{array}\right.
\end{align*}
The same result holds for $C_k^+$ except the case $(k,r)=(4,3)$, in which
\begin{align*}ex(n,3,C_4^+)={{n}\choose {r}}-{{n-1}\choose r}+\max\{n-3,4\lfloor\frac{n-1}{4}\rfloor\}.
\end{align*}
\end{lemma}

In \cite{Furedi&Jiang&Seiver2014}, F\"{u}redi, Jiang and Seiver determined $ex(n,r,\mathcal{P}_k)$ for $r\geq 3$ and sufficiently large $n$. F\"{u}redi and Jiang \cite{Furedi&Jiang2014} determined $ex(n,r,\mathcal{C}_k)$ for all $k\geq 3$, $r\geq 4$ and sufficiently large $n$. Kostochka, Mubayi and Verstra\"{e}te \cite{Kostochka&Mubayi&Verstraete2015} extended their results to all $r\geq 3$.
\begin{lemma}[\cite{Furedi&Jiang2014,Furedi&Jiang&Seiver2014,Kostochka&Mubayi&Verstraete2015}]\label{turanloosePkCk} Let $\delta,$ $\varepsilon>0$ and $r, k, t, n$ be positive integers such that $r\geq 3,$ $t=\lfloor \frac{k-1}{2}\rfloor$ and $1/n\ll \delta\ll \varepsilon$. Then
\begin{align*}ex(n,r,\mathcal{P}_k)={{n}\choose {r}}-{{n-t}\choose r}+\left\{\begin{array}{ll}
0, & \mbox{if $k$ is odd,}\\ 1, & \mbox{if $k$ is even}.
\end{array}\right.
\end{align*}
And
\begin{align*}ex(n,r,\mathcal{C}_k)={{n}\choose {r}}-{{n-t}\choose r}+\left\{\begin{array}{ll}
0, & \mbox{if $k$ is odd,}\\ 1, & \mbox{if $k$ is even and $k\neq 4$}.
\end{array}\right.
\end{align*}
Especially, for $k=4$,
\begin{align*}ex(n,r,\mathcal{C}_4)={{n}\choose {r}}-{{n-1}\choose r}+\big\lfloor\frac{n-1}{s}\big\rfloor.
\end{align*}
\end{lemma}

For an $r$-graph $H$, let $\partial H$ denote the \emph{shadow} of $H$, which is the $(r-1)$-graph consisting of $(r-1)$-sets contained in some edge of $H$.

\begin{lemma}[\cite{Kostochka&Mubayi&Verstraete2015}]\label{stability}
For fixed $r\geq 3$ and $k\geq4$, let $t=\lfloor\frac{k-1}{2}\rfloor$ and $H$ be an $r$-graph on $n$ vertices with $|E(H)|\sim t{{n}\choose{r-1}}$ containing no $P_k^+$ or $C_k^+$. Then there exists $G\subset \partial H$ with $|E(G)|\sim {{n}\choose{r-1}}$ and a set $L$ of $t$ vertices of $H$ such that $L\cap V(G)=\emptyset$ and $e\cup \{v\}\in E(H)$ for any $(r-1)$-edge $e\in G$ and any $v\in L$.
\end{lemma}

Suppose $H$ is an $r$-graph which satisfies the conditions of Lemma \ref{stability}. Let $L=\{v_1,\ldots,v_{t}\}$, $\widetilde{H}=H-L$. Let $\tau$ be a given positive integer. For two vertices $u, v\in V(H)$, denote the number of edges in $H$ containing $u$ and $v$ by $d_H(u,v)$. A pair of vertices $\{u,v\}$ with $u\in L$ and $v\notin L$ is $\tau$-\emph{small} if $d_H(u,v)\leq \tau {n\choose {r-3}}$, and  \emph{$\tau$-big} otherwise. For $v\in V(\widetilde{H})$, define the $\tau$-\emph{small degree} of $v$, denoted by $d_s(v,L,\tau)$, to be the number of vertices in $L$ which could form a $\tau$-small pair with $v$, i.e., $d_s(v,L,\tau)=|\{i\in [t]: \{v_i,v\} \mbox{ is $\tau$-small}\}|$. An edge $e\in E(H)$ is \emph{crossing} if $|e\cap L|=1$, and denoted by ${\rm Cross}(H)$ the set of all crossing edges in $H$. An $r$-set of vertices $e$ is  a \emph{missing edge} of $H$ if $e\notin E(H)$ and $|e\cap L|=1$, and let ${\rm Miss}(H)$ denote the set of all missing edges of $H$. Notice that all these definitions depend on $L$.

For $F=\{e_1,\ldots,e_k\}\in \mathcal{P}_k$, the \emph{end edges} of $F$  are the two edges, say $e_1$ and $e_k$, that intersect only one edge in $F$. A vertex $v$ is an \emph{end point} of $F$ if $v$ is in an end edge of $F$ and $d_{F}(v)=1$, where $d_{F}(v)$ is the number of edges in $F$ containing $v$. Denote by ${\rm end}(F)$ the set of all the end points of $F$. A pair of vertices $\{u,v\}$ is  an \emph{end pair} of  $F$ if $u,v\in {\rm end}(F)$ and $u,v$ belong to different end edges of $F$. Let $D(F)$ denote the set of all the end pairs of $F$. For any subgraph $F$ of $H$, let $C(F)$ denote the set of colors on the edges of $F$.

The following lemma can help us to reduce the problem from finding long loose paths or cycles to finding short loose paths, due to the degrees of some vertex pairs are big enough.

\begin{lemma}\label{findckpk} Let $r,\ell,t,\tau$ be integers such that $r\geq 3$, $\ell\geq 1$, $t\geq 1$ and $\tau=r(\ell+2t)$. Given an edge-colored $K_n^r$  and a rainbow subgraph $H$ of $K_n^r$. Let    $L=\{v_1,\ldots,v_t\}$ be a set of $t$ vertices in $H$. Let $S=\{v\in V(\widetilde{H}): d_s(v,L,\tau)\geq 1\}$, $\bar{S}=V(H)\backslash (L\cup S)$. Suppose that $|S|\leq n/2$ and  there exists a rainbow loose path of length $\ell$, say $F$, in $\widetilde{H}$. Then

${\rm (i)}$  if there is an end pair $\{x,y\}\in D(F)$ such that $x,y\in \bar{S}$, then $K_n^r$ contains all rainbow loose cycles of length $\ell+2i$, $i\in [t]$;

${\rm (ii)}$  if there is an end point $x\in {\rm end}(F)$ such that $x\in \bar{S}$, then $K_n^r$ contains all rainbow loose paths of length $\ell+2i$, $i\in [t]$.

\end{lemma}

\begin{prof}[Proof of Lemma 2.7.] We can see that for a $\tau$-big pair $\{u,v\}$, any color set $C$ with $|C|\leq \ell+2t$ and any $U\subset V(K_n^r)\setminus\{u,v\}$ with $|U|\leq \tau$, there always exists an edge $e\in E(H)$ such that $\{u,v\}\subset e$, $e\cap U=\emptyset$ and $C(e)\notin C$. Since $|S|\leq n/2$, we can choose $u_1,u_2,\ldots,u_{t}\in \bar{S}\backslash V(F)$.

 ${\rm (i)}$  Let $W=\{x,y,u_1,\ldots,u_{t}\}$. By the definition of $S$, all pairs in $\{\{v,u\}:v\in L,u\in W\}$ are $\tau$-big.  For any $i\in [t]$,  consider the $2i$ $\tau$-big pairs $\{x,v_1\}$, $\{v_1,u_1\}$, $\{u_1,v_2\}$, $\{v_2,u_2\}$,\ldots, $\{v_{i-1},u_{i-1}\}$, $\{u_{i-1},v_i\}$,  $\{v_{i},y\}$. Then there is a rainbow loose path of length $2i$ in $H$, say $F_1$, with $x,y $ being one of the end pairs and $V(F)\cap V(F_1)=\{x, y\}$. Thus $F\cup F_1$ forms a rainbow loose cycle of length $\ell+2i$, as desired.

 ${\rm (ii)}$   Let $W=\{x,u_1,\ldots,u_{t}\}$. By the definition of $S$, all pairs in $\{\{v,u\}:v\in L,u\in W\}$ are $\tau$-big.  For any $i\in [t]$,  consider the $2i$ $\tau$-big pairs $\{x,v_1\}$, $\{v_1,u_1\}$, $\{u_1,v_2\}$, $\{v_2,u_2\}$,\ldots, $\{v_{i-1},u_{i-1}\}$, $\{u_{i-1},v_i\}$, $\{v_{i},u_i\}$. Then there is a rainbow loose path of length $2i$ in $H$, say $F_2$,  with $x$ being one of the end points and $V(F)\cap V(F_2)=\{x\}$. Thus $F\cup F_2$ forms a rainbow loose path of length $\ell+2i$, as desired.
\end{prof}

In the following, for a fixed family of $r$-graphs $\mathcal{F}$, we say $H$ is $\mathcal{F}$-\emph{free} if $H$ contains no hypergraph in $\mathcal{F}$ as a subgraph. Given an edge-coloring of $H$, $H$ is \emph{rainbow $\mathcal{F}$-free} if $H$ contains  no rainbow hypergraph in $\mathcal{F}$ as a subgraph.

\section{even paths or cycles}

For the lower bound, consider the following edge-coloring of $K_{n}^{r}$ with ${n\choose r}-{{n-t+1}\choose r}+1$ colors. Take a set $L$ of $t-1$ vertices and color all edges containing some vertex in $L$ with different colors, which uses ${n\choose r}-{{n-t+1}\choose r}$ colors. Then color all the remaining edges with one additional color. We can see that such an edge-coloring yields no rainbow loose cycles or loose paths of length $k$ with $k=2t$. So $ar(n,r,\mathcal{C}_{k})\geq {n\choose r}-{{n-t+1}\choose r}+2$ and $ar(n,r,\mathcal{P}_{k})\geq {n\choose r}-{{n-t+1}\choose r}+2$.

For the upper bound, we have $ar(n,r,\mathcal{C}_{k})\leq ar(n,r,C_k^{+})={n\choose r}-{{n-t+1}\choose r}+2$ and $ar(n,r,\mathcal{P}_{k})\leq ar(n,r,P_k^{+})={n\choose r}-{{n-t+1}\choose r}+2$ by Lemma \ref{linear}, since a linear path or cycle is also a loose path or cycle.

This completes the proof.

\section{odd paths or cycles}

For the lower bound, consider the following edge-coloring of $K_{n}^{r}$ with ${n\choose r}-{{n-t+1}\choose r}+2$ colors. Take a set $L$ of $t-1$ vertices and color all edges containing some vertex in $L$ with different colors, which uses ${n\choose r}-{{n-t+1}\choose r}$ colors. Then color all the remaining edges with two additional colors. We can see that such an edge-coloring yields no rainbow loose cycles or loose paths of length $k$ with $k=2t+1$. So $ar(n,r,\mathcal{C}_{k})\geq {n\choose r}-{{n-t+1}\choose r}+3$ and $ar(n,r,\mathcal{P}_{k})\geq {n\choose r}-{{n-t+1}\choose r}+3$.

For the upper bound, we will first find a subgraph in $K_{n}^{r}$ satisfying the conditions of Lemma \ref{stability} to get a nice structure, then analyse the quantitative relationships corresponding to the structure and finally try to find long loose paths or cycles to get contradictions.

Choose positive constants $\varepsilon_1, \varepsilon_2, \varepsilon_3, \varepsilon_4, \varepsilon_5$ satisfying $1/n\ll \varepsilon_5\ll \varepsilon_4\ll \varepsilon_3\ll \varepsilon_2\ll \varepsilon_1\ll k, r$. Suppose that there is an edge-coloring of $K_{n}^{r}$ with ${n\choose r}-{{n-t+1}\choose r}+3$ colors yielding no rainbow loose paths or loose cycles of length $k$. Let $H$ be the rainbow subgraph of $K_{n}^{r}$ obtained by choosing one edge from each color class. Obviously, $|E(H)|={n\choose r}-{{n-t+1}\choose r}+3$. Notice that $H$ is clearly $\mathcal{P}_{k}$-free and $\mathcal{C}_{k}$-free. Next, we will find a subgraph $H'$ of $H$ containing no copy of $P_{k-1}^{+}$ to apply Lemma \ref{stability}.

Let $\mathbf{\mathcal{P}}$ be any maximal collection of pairwise edge-disjoint copies of $P_{k-1}^+$ in $H$. Then if we delete all the edges in $\mathbf{\mathcal{P}}$, the remaining hypergraph is clearly $P_{k-1}^+$-free.  We will first claim that the size of $\mathbf{\mathcal{P}}$ is at most $O(n)$ and therefore the total number of deleted edges of $H$ is at most $O(n)$.

\begin{claim}[\cite{Tang&Li&Yan2022}]\label{numberP_{k-1}^{+}}
\quad $|\mathbf{\mathcal{P}}|\leq \frac{k'-1}{2} n$, where $k'=4((k-2)(r-1)+2)$.
\end{claim}

Now let $H'$ be the subgraph by deleting from $H$ all the edges in each $F\in \mathbf{\mathcal{P}}$. We have $H'$ is $P_{k-1}^+$-free due to $\mathbf{\mathcal{P}}$ is maximal, and
\begin{align*}
|E(H')|&\geq |E(H)|-\bigg(\frac{k'-1}{2} \bigg)n(k-1)\\
 &\geq (t-1){n\choose{r-1}}-\delta n^{r-1},
\end{align*}
where the last inequality is due to $r\geq 3$ and $1/n\ll \delta$.

By Lemma \ref{stability}, there exists a set of $t-1$ vertices $L$ such that $d_{H'}(v,V(H')\backslash L)\geq {n-1\choose {r-1}}-\varepsilon_5 n^{r-1}$. Then $d_H(v,V(H)\backslash L)\geq {n-1\choose {r-1}}-\varepsilon_5 n^{r-1}$. Now we will analyse the structure of $H$ in detail to obtain a contradiction. Let $\tau= r(2t+1)$. Recall that $\widetilde{H}=H-L$, $S=\{v\in V(\widetilde{H}): d_s(v,L,\tau)\geq 1\}$ and $\bar{S}=V(H)\backslash (L\cup S)$.

\noindent{\bf Claim 4.2.} \quad The number of $\tau$-small pairs is at most $\varepsilon_3 n$. In particular, $|S|\leq \varepsilon_3 n$.

\begin{prof}[Proof of Claim 4.2.] Let $s$ be the number of $\tau$-small pairs. We prove this claim by a simple double counting on $|{\rm Miss}(H)|$. Notice that $|{\rm Miss}(H)|+|{\rm Cross}(H)|=(t-1){{n-t+1}\choose {r-1}}$. According to Lemma \ref{stability}, $|{\rm Cross}(H)|\geq (t-1)\big[{{n-1}\choose {r-1}}-\varepsilon_5 n^{r-1}\big]$, and so
\begin{align}\label{1}
|{\rm Miss}(H)|\leq \varepsilon_4 n^{r-1}.
\end{align}

On the other hand, one $\tau$-small pair corresponds to at least ${{n-t}\choose {r-2}}-\tau{n\choose {r-3}}$ missing edges. Each missing edge contains at most $(r-1)$ $\tau$-small pairs, thus
\begin{align}\label{2}
|{\rm Miss}(H)|\geq \frac{\big({{n-t}\choose {r-2}}-\tau{n\choose {r-3}}\big) s}{r-1}.
\end{align}
Therefore, by (\ref{1}) and (\ref{2}), $s\leq \varepsilon_3 n$. In particular, since each vertex in $S$ contributes to at least one $\tau$-small pair, we have $|S|\leq s\leq \varepsilon_3 n$.
\end{prof}

Next, we partition the edges in $\widetilde{H}$. Let $E_{i}=\{e\in\widetilde{H}: |e\cap \bar{S}|=i\}$ and $E_{[i,j]}=\{e\in\widetilde{H}: i\leq|e\cap \bar{S}|\leq j\}$. Then $E(\widetilde{H})=\bigcup_{i=0}^{r}E_{i}=E_{[0,r]}$. We will show that $|E_0|$, $|E_{[2,r-1]}|$ and $|E_r|$ are relatively small with respect to $|E_1|$.

\noindent{\bf Claim 4.3.} \quad There exists a positive integer $M$, such that $|E_{0}|\leq \max\{{{|S|}\choose {r}}-{{|S|-t}\choose r},{M\choose r}\}$.

\begin{prof}[Proof of Claim 4.3.] Since $K_n^r$ is rainbow $\mathcal{C}_k$-free or rainbow $\mathcal{P}_k$-free, we have $E_0$ is  $\mathcal{C}_k$-free or $\mathcal{P}_k$-free. This, together with Lemma \ref{turanloosePkCk}, implies the desired bound.
\end{prof}

\noindent{\bf Claim 4.4.} \quad $|E_{r}|\leq 2$.

\begin{prof}[Proof of Claim 4.4.] Since $K_n^r$ is rainbow $\mathcal{C}_k$-free or rainbow $\mathcal{P}_k$-free and $|S|\leq \varepsilon_3 n$, we have $E_r$ is  $\mathcal{P}_3$-free by Lemma \ref{findckpk}. This, together with Lemma \ref{arrP3}, implies the desired bound.
\end{prof}

\noindent{\bf Claim 4.5.} \quad $|E_{i}|\leq 2{{|S|}\choose {r-i}}$, for every $i\in \{2,3,\ldots, r-1\}$.

\begin{prof}[Proof of Claim 4.5.] We only need to show that for every $(r-i)$-tuple in ${{S}\choose {r-i}}$, there are at most two edges in $E_{i}$ containing this tuple. Suppose to the contrary that there exist three edges $e_1, e_2, e_3$ in $E_i$ containing some common $(r-i)$-tuple in ${{S}\choose {r-i}}$, for some $i\in \{2,3,\ldots, r-1\}$.

\noindent {\bf Case 1}. $(e_1\cap \bar{S})\setminus (e_2\cup e_3)\neq \emptyset$, $(e_2\cap \bar{S})\setminus (e_1\cup e_3)\neq \emptyset$ and $(e_3\cap \bar{S})\setminus (e_1\cup e_2)\neq \emptyset$.

Choose $u_1\in (e_1\cap \bar{S})\setminus (e_2\cup e_3)$, $u_2\in (e_2\cap \bar{S})\setminus (e_1\cup e_3)$, $u_3\in (e_3\cap \bar{S})\setminus (e_1\cup e_2)$. Then take $e_1'=\{u_1, w', w_1, w_2,\ldots, w_{r-2}\}$, where $w', w_1, w_2,\ldots, w_{r-2}\in \bar{S}\setminus (e_1\cup e_2\cup e_3)$, $e_2'=\{u_2, w', w'', v_1, v_2,\ldots, v_{r-3}\}$, where $w'', v_1, v_2,\ldots, v_{r-3}\in \bar{S}\setminus (e_1\cup e_2\cup e_3\cup e_1')$ and $e_3'=\{u_3, w'', x_1, x_2,\ldots, x_{r-2}\}$, where $x_1, x_2,\ldots, x_{r-2}\in \bar{S}\setminus (e_1\cup e_2\cup e_3\cup e_1'\cup e_2')$. Since $e_1, e_2, e_3\in E_i$, $e_1, e_2, e_3$ have different colors. We can see that $\{e_1', e_1, e_2\}$ and $\{e_1', e_1, e_3\}$ are two loose paths of length 3, so they cannot be rainbow by Lemma \ref{findckpk}. Thus, $C(e_1')=C(e_1)$. Similarly, we have $C(e_2')=C(e_2)$ according to that $\{e_2', e_2, e_3\}$ and $\{e_2', e_2, e_1\}$ are two loose paths of length 3, and $C(e_3')=C(e_3)$ according to that $\{e_3', e_3, e_2\}$ and $\{e_3', e_3, e_1\}$ are two loose paths of length 3. But $\{e_1', e_2', e_3'\}$ is also a loose path of length 3 and it is rainbow. Therefore, $K_{n}^r$ contains a rainbow loose path and a loose cycle of length $k$ by Lemma \ref{findckpk}, a contradiction.

\noindent {\bf Case 2}. $(e_1\cap \bar{S})\setminus (e_2\cup e_3)= \emptyset$.

Choose $u_1\in e_1\cap e_2\cap \bar{S}\setminus e_3$, $u_2\in e_1\cap e_3\cap \bar{S}\setminus e_2$. Then take $e_1'=\{u_1, w_1, w_2,\ldots, w_{r-1}\}$, where $w_1, w_2,\ldots, w_{r-1}\in \bar{S}\setminus (e_1\cup e_2\cup e_3)$ and $e_2'=\{u_2, v_1, v_2,\ldots, v_{r-1}\}$, where $v_1, v_2,\ldots, v_{r-1}\in \bar{S}\setminus (e_1\cup e_2\cup e_3\cup e_1')$. Since $e_1, e_2, e_3\in E_i$, $e_1, e_2, e_3$ have different colors. We can see that $\{e_1', e_1, e_3\}$ and $\{e_1', e_2, e_3\}$ are two loose paths of length 3, so they cannot be rainbow by Lemma \ref{findckpk}. Thus $C(e_1')=C(e_3)$. Similarly, we have $C(e_2')=C(e_2)$ according to that $\{e_2', e_1, e_2\}$ and $\{e_2', e_3, e_2\}$ are two loose paths of length 3. But $\{e_1', e_1, e_2'\}$ is also a loose path of length 3 and it is rainbow. Therefore, $K_{n}^r$ contains a rainbow loose path and a loose cycle of length $k$ by Lemma \ref{findckpk}, a contradiction.
\end{prof}

\noindent{\bf Claim 4.6.} \quad $|E(\widetilde{H})|\geq \frac{|S|}{r-1}\big[{{n-t}\choose {r-2}}-\tau{n\choose {r-3}}\big]+3$.

\begin{prof}[Proof of Claim 4.6.] Since $|S|\leq s$, we only need to prove $|E(\widetilde{H})|\geq |{\rm Miss}(H)|+3$ by (\ref{2}).

Define $E(L)=\{e\in {{V(K_n^r)}\choose r}: |e\cap L|\geq2\}$, $E_H(L)=E(L)\cap E(H)$. Notice that
\begin{align}
& |E(H)|=|E(\widetilde{H})|+|{\rm Cross}(H)|+|E_H(L)|,\label{a}\\
&|{\rm Cross}(H)|+|{\rm Miss}(H)|={n\choose r}-{{n-t+1}\choose r}-|E(L)|,\label{b}\\
&|E(H)|\geq {n\choose r}-{{n-t+1}\choose r}+3. \label{c}
\end{align}
Combining (\ref{a}), (\ref{b}) and (\ref{c}), we get
\begin{align*}
|E(\widetilde{H})|&\overset{(\ref{a}), (\ref{c})}{\geq} {n\choose r}-{{n-t+1}\choose r}+3-|{\rm Cross}(H)|-|E_H(L)|\\
&\overset{(\ref{b})}{=}|{\rm Cross}(H)|+|{\rm Miss}(H)|+|E(L)|-|E_H(L)|-|{\rm Cross}(H)|+3\\
&\geq |{\rm Miss}(H)|+3.
\end{align*}
\end{prof}

It is easy to see that $|S|\geq 1$. Otherwise $\bar{S}=V\backslash L$. By Claim 4.6, $|E(\widetilde{H})|\geq 3$. By Lemma \ref{arrP3}, there is a rainbow loose path of length 3 in the subgraph of edge-colored $K_n^r$ induced by $V\backslash L$, which implies that there is a rainbow loose path of length $k$ and a loose cycle of length $k$ in $K_n^r$ by Lemma \ref{findckpk}, a contradiction.

By Claims 4.3, 4.4, 4.5 and 4.6, we have
\begin{align}
|E_1|&=|E(\tilde{H})|-|E_0|-|E_{[2,r-1]}|-|E_r|\notag\\
&\geq \frac{|S|}{r-1}\bigg[{{n-t}\choose {r-2}}-\tau{n\choose {r-3}}\bigg]+3-|E_0|-|E_{[2,r-1]}|-|E_r|\notag\\
&\geq \frac{|S|}{r-1}\bigg[{{n-t}\choose {r-2}}-\tau{n\choose {r-3}}\bigg]-(t+1){{|S|}\choose {r-1}} \notag\\
& \geq \frac{|S|}{r-1}\bigg[{{n-t}\choose {r-2}}-\tau{n\choose {r-3}}\bigg]-|S|\frac{t+1}{(r-1)!}\varepsilon_3^{r-2}n^{r-2}\notag\\
&\geq \frac{|S|}{r-1}\bigg[{{n-t}\choose {r-2}}-\tau{n\choose {r-3}}-\frac{t+1}{(r-2)!}\varepsilon_3^{r-2}n^{r-2}\bigg]\notag\\
& \geq \frac{|S|}{r-1}\bigg[{{n-t}\choose {r-2}}-\varepsilon_2 n^{r-2}\bigg], \label{*}
\end{align}
where the third  inequality follows from $|S|\leq \varepsilon_3 n$ by Claim 4.2.

Next, we partition the shadow of $E_1$ within $S$ into two parts, i.e., $\partial(E_1)\cap {S\choose {r-1}}=F_1\cup F_{2^+}$, where $F_1=\{f\in \partial(E_1)\cap {S\choose {r-1}}: d_{E_1}(f)=1 \}$, $F_{2^+}=\{f\in \partial(E_1)\cap {S\choose {r-1}}: d_{E_1}(f)\geq 2 \}$ and $d_{E_1}(f)=|\{e\in E_1: f\subset e\}|$.

We will show that $|F_{2^+}|$ is relatively small, which deduces that the size of $F_1$ is large enough to lead to a contradiction.

\noindent {\bf Claim 4.7.}  $|F_{2^+}|\leq  ex(|S|,r-1,\mathcal{P}_2)$.

\begin{prof}[Proof of Claim 4.7.] \quad Suppose $|F_{2^+}|> ex(|S|,r-1,\mathcal{P}_2)$. Then there is a loose path of length 2 in $F_{2^+}$, say $\{f_1,f_2\}$. Let $u_1=f_1\setminus f_2$. By the definition of $F_{2^+}$,  there exist $x,x',y,y'\in \bar{S}$  such that  $\{x\}\cup f_1$, $\{y\}\cup f_1$, $\{x'\}\cup f_2$ and $\{y'\}\cup f_2$ are four edges in $E_1$. Choose $e=\{z_1,\ldots,z_{r-1},u_1\}$, where $z_i\in \bar{S}\backslash \{x,y,x',y'\}$, $i\in [r-1]$. Since $E_1\subset E(\widetilde{H})$, the four edges $\{x\}\cup f_1$, $\{y\}\cup f_1$, $\{x'\}\cup f_2$ and $\{y'\}\cup f_2$ have different colors and at least three of them have different colors from $C(e)$. Therefore, we can take two edges among these three edges containing $f_1$ and $f_2$  respectively, together with $e$, to form a rainbow loose path of length 3 with an end pair in $\bar{S}$. Thus, $K_{n}^r$ contains a rainbow loose path of length $k$ and a loose cycle of length $k$ by Lemma \ref{findckpk}, a contradiction.
\end{prof}

\noindent {\bf Claim 4.8.} $|F_1|\geq |S|\varepsilon_1 n^{r-2}$.

\begin{prof}[Proof of Claim 4.8.] \quad
Notice that $|E_1|=\sum_{f\in F_1}d_{E_1}(f)+\sum_{f\in F_{2^+}}d_{E_1}(f)\leq |F_1|+n|F_{2^+}|$. This, combined with Claim 4.7 and (\ref{*}), implies that
\begin{align*}
|F_1|&\geq \frac{|S|}{r-1}\bigg[{{n-t}\choose {r-2}}-\varepsilon_2 n^{r-2}\bigg]-ex(|S|,r-1,\mathcal{P}_2)\cdot n.\\
\end{align*}
Obviously, $ex(|S|,r-1,\mathcal{P}_2)\leq \frac{|S|}{r-1}$. Then
 \begin{align*}
|F_1|&\geq \frac{|S|}{r-1}\bigg[{{n-t}\choose {r-2}}-\varepsilon_2 n^{r-2}\bigg]-\frac{|S|}{r-1}\cdot n\\
&\geq \frac{|S|}{r-1}\bigg[{{n-t}\choose {r-2}}-\varepsilon_2 n^{r-2}-n]\\
&\geq |S|\varepsilon_1 n^{r-2}.
\end{align*}
\end{prof}

However, on the other hand, $|F_1|\leq {{|S|}\choose {r-1}}$. It follows that
\begin{align*}
&|S|\varepsilon_1 n^{r-2}\leq |F_1|\leq {{|S|}\choose {r-1}}\leq \frac{|S|^{r-1}}{(r-1)!}\\
&\Rightarrow\quad \varepsilon_1 n^{r-2}\leq \frac{|S|^{r-2}}{(r-1)!}\\
&\Rightarrow\quad \varepsilon_2n <|S|,
\end{align*}
which contradicts to $|S|\leq \varepsilon_3 n $ by Claim 4.2.\\
This completes the proof.

\section*{References}

\bibliography{bibl}

\end{document}